\documentclass[11pt,a4paper]{article}

\usepackage{amsmath}
\usepackage{amsfonts}

\usepackage[latin1]{inputenc}
\usepackage{fontenc}
\usepackage{euscript}
\usepackage{array}
\usepackage{vmargin}
\usepackage{amssymb}
\usepackage{latexsym}
\usepackage{enumerate}
\usepackage{multirow}
\usepackage{graphicx}

\newcommand{\pg}{\hspace{0.6cm}}

\def\N{{\rm I\kern-.20em N}}
\def\R{{\rm I\kern-.20em R}}
\def\indi{{1\kern-.20em\rm I}}

\newtheorem{definition}{Definition}[section]
\newtheorem{proposition}{Proposition}[section]
\newtheorem{remark}{Remark}

\newcommand {\edem}{\hfill $\square$ \end {proof}}
\newtheorem{ex}{Example}[section]
\newcommand{\bdem} {\begin{proof}}

\begin{document}
\bibliographystyle{plain}
\hyphenation{o-pe-ra-tors}

\title{Stability and contagion measures for spatial extreme value analyses}
\author{Fonseca, C.\footnote{\noindent Instituto Politécnico da Guarda, Portugal, E-mail: {\tt
cfonseca@ipg.pt}\hspace*{2mm}},  Ferreira, H., Pereira, L. and
Martins, A.P. \footnote{\noindent Departamento de Matemática,
Universidade da Beira Interior, Portugal, E-mail: {\tt helena.ferreira@ubi.pt}, {\tt
lpereira@ubi.pt}, {\tt
amartins@ubi.pt}\hspace*{2mm}}}


\date{}
\maketitle

\noindent {\bf Abstract:} As part of global climate change an accelerated hydrologic cycle (including an increase in heavy precipitation) is anticipated (\cite{trenb1},\cite{trenb2}). So, it is of great importance to be able to quantify high-impact hydrologic relationships, for example, the impact that an extreme precipitation (or temperature) in a location has on a surrounding region. Building on the Multivariate Extreme Value Theory we propose a contagion index and a stability index. The contagion index makes it possible to quantify the effect that an exceedance above a high threshold can have on a region. The stability index reflects the expected number of crossings of a high threshold in a region associated to a specific location ${\bf i}$, given the occurrence of at least one crossing at that location. We will find some relations with well-known extremal dependence measures found in the literature, which will provide immediate estimators. For these estimators an application to the annual maxima precipitation in Portuguese regions is presented.





\pagenumbering{arabic}

\section{Introduction}

\pg The need to model and predict environmental extreme events such as hurricanes, floods, droughts, heat waves and other high impact events, which can lead to a devasting impacts, ranging from disturbances in ecosystems to economic impacts on society as well as loss of life, motivated the modeling of spatial extremes.

A common method of modeling spatial extremes is through max-stable processes, which are an infinite dimensional generalization of multivariate extreme value distributions (Haan (1984) \cite{haan1}, Vatan (1985) \cite{vat}, de Haan and Pickands (1986) \cite{haan2}).

Max-stable processes can be, for example, good approximations for annual maxima of daily spatial rainfall (Smith \cite{smi}, Coles \cite{col}, Schlather \cite{sch1}, among others) and therefore have been widely applied to real data.

Briefly, a max-stable process ${\bf X}=\left\{X_{\mathbf{i}}\right\}_{{\mathbf{i}}\in \mathbb R^d}$ is the limit process of maxima of i.i.d. random fields $\left\{Y^{(j)}_{\mathbf{i}}\right\}_{{\mathbf{i}}\in \mathbb R^d},$ $j\geq 1.$ Namely, for suitable $\left\{a_n({\mathbf{i}})>0\right\}_{n\geq 1}$ and $\left\{b_n({\mathbf{i}})\right\}_{n\geq 1}$ sequences of real constants, $$X_{{\mathbf{i}}}=\lim_{n\to \infty}\frac{\bigvee_{j=1}^nY^{(j)}_{\mathbf{i}}-b_n({\mathbf{i}})}{a_n({\mathbf{i}})},\qquad {\mathbf{i}}\in \R^d,$$ provided the limit exists.

We shall consider $d=2,$ that is ${\bf X}=\left\{X_{\mathbf{i}}\right\}_{{\mathbf{i}}\in \mathbb R^2}.$ The distribution of $(X_{\mathbf{i}_1},\ldots,X_{\mathbf{i}_k})$ is a Multivariate Extreme Value (MEV) distribution $G_{{\bf A}}$, ${\bf{A}}=\{{\bf{i}}_1,\ldots,{\bf{i}}_k\}$, and we can assume, without loss of generality, that the margins of ${\bf X}$ have a unit Fréchet distribution, $F(x)=\exp(-x^{-1}),\ x>0$ (Resnick \cite{res}). The distribution $G_{{\bf A}}$ can then be defined by
\begin{equation}
G_{{\bf A}}(x_1,\ldots,x_k)=\exp(-V_{\bf A}(x_1,\ldots,x_k)),\quad x_i\in \R^+,\quad i=1,\ldots,k,
\end{equation}
where $V_{{\bf A}}$ denotes the exponent function of the MEV distribution $G_{{\bf A}}$.

The exponent function summarizes the extremal dependence structure of $G_{{\bf A}}$ and the scalar $V_{\bf{A}}(1,\ldots,1)$  defines the extremal coefficient $\epsilon_{\bf{A}}$ detailed in Schlather and Tawn (2003,\cite{sch}), which summarizes the extremal dependence between the variables indexed in the region ${\bf{A}}$. This coefficient takes values between 1 and $k$, with a value of 1 corresponding to complete dependence and a value of $k$ corresponding to complete independence. Its value can be thought of as the number of effectively independent locations among the $k$ under consideration.


If we consider ${\bf A}=\{{\bf{i}},{\bf{j}}\}$ we find the extremal coefficient of Tiago de Oliveira (1962/63, \cite{tiago}) which is related with the bivariate upper tail dependence coefficient $\lambda_{\{\bf {i},\bf{j}\}}=\lim_{u\uparrow 1}P\left(F(X_{\bf{i}})>u|F(X_{\bf{j}})>u\right)$, introduced in Sibuya (1960, \cite{sibuya}), as $\epsilon_{\{\bf {i},\bf{j}\}}=2-\lambda_{\{\bf {i},\bf{j}\}}$.

Although these measures are very useful to analyze the dependence among extremal events, there remain important questions to be answered, for example, the influence of an extreme event on the regional smoothness of $\bf X$ and the contagion effect of an extreme event at a specific location over the variables of $\bf X$ indexed in a region of $\mathbb{R}^2$.


In Section 2 we propose a measure, called contagion index, that enables to quantify the impact that an exceedance of a high threshold can have on a region and we present its relation with bivariate extremal coefficients.

Clearly an extreme event could affect the smoothness of a random field over a region so, in Section 3 we propose a stability index on a region $\mathbf{A}$ associated to a specific location $\mathbf{i}$, defined as the expected number of crossings of a high threshold u in $\mathbf{A}$ associated with $\mathbf{i}$, given that there is at least one crossing in $\mathbf{A}$ associated to $\mathbf{i}$. We also present some properties of this coefficient.

Section 4 is devoted to illustrate the previous measures in a max-stable M4 random field.

Based on relations of our indices with well-known dependence measures, for which estimators and respective properties have already been studied in the literature, in Section 5 we present estimators for the stability and contagion indices. We end with an application to the annual maxima precipitation in Portuguese regions.

\section{Contagion Index}

\pg The occurrence of an extreme event at a given location $\mathbf{i}$ may spread throughout a region of locations. In this section we define a measure for assessing the effect of an exceedance above a high threshold $u$ at a specific location $\mathbf{i}$ on a region $\mathbf{A}$ of locations.

 \begin{definition}
 Let $X = \{X_{\bf{i}}\}_{{\bf{i}}\in\R^2}$ be a max-stable random field with unit Fréchet margins $F$ and $\mathbf{A}$ a region of $\mathbb{R}^2$. The contagion index from the location $\mathbf{i}$ to the region $\mathbf{A}$ is defined as
\begin{equation}
CI(\mathbf{A},\mathbf{i})=\lim_{u\uparrow 1} E\left(\sum_{\mathbf{j}\in \mathbf{A}}\indi_{\{F(X_{\mathbf j})>u\}} \Biggl|\ F(X_{\mathbf i})>u\right). \label{IC}
\end{equation}
\end{definition}\vspace{0.3cm}

The $CI(\mathbf{A},\mathbf{i})$ is the conditional expected number of exceedances above a high threshold $u$ in $\mathbf{A}$, given $X_{\mathbf i}$ exceeds $u$, that is, the $CI(\mathbf{A},\mathbf{i})$ measures the impact that the event $\left\{X_{\mathbf i}>u\right\}$ has on the region $\mathbf{A}$.

We remark that the conditioning location ${\bf i}$ does not necessarily have to be in the region ${\bf A}$.

The following proposition states that $CI(\mathbf{A},\mathbf{i})$ relates with the bivariate extremal dependence \linebreak coefficients $\epsilon_{\{\mathbf{i},\mathbf{j}\}}$, ${\mathbf{j}}\in \mathbf{A}$, and the tail dependence coefficients $\lambda_{\{\mathbf{i},\mathbf{j}\}}$, ${\mathbf{j}}\in \mathbf{A}$.

\begin{proposition}
For any max-stable random field with unit Fréchet margins $F$, ${\mathbf i}\in \mathbb{R}^2$ and ${\mathbf{A}}\subset \mathbb{R}^2$, we have
$$
CI(\mathbf{A},\mathbf{i})=\sum_{\mathbf{j}\in \mathbf{A}}\lambda_{\{{\mathbf {i}},{\mathbf{j}}\}}=2\left|{\mathbf{A}}\right|-\sum_{\mathbf{j}\in \mathbf{A}}\epsilon_{\{{\mathbf {i}},{\mathbf{j}}\}}.
$$
\end{proposition}

\bdem Observe that
\begin{eqnarray*}
CI(\mathbf{A},\mathbf{i})&=&\sum_{\mathbf{j}\in \mathbf{A}}\lim_{u\uparrow 1} P\left(F(X_{\mathbf {j}})>u\left| F(X_{\mathbf {i}})>u\right.\right)\\[0.3cm]
&=&\sum_{\mathbf{j}\in \mathbf{A}}\lambda_{\{{\mathbf {i}},{\mathbf{j}}\}}= \sum_{\mathbf{j}\in \mathbf{A}}\left(2-\epsilon_{\{{\mathbf {i}},{\mathbf{j}}\}}\right)= 2\left|{\mathbf{A}}\right|-\sum_{\mathbf{j}\in \mathbf{A}}\epsilon_{\{{\mathbf {i}},{\mathbf{j}}\}}.
\end{eqnarray*}
\edem\vspace{0.3cm}

An $CI(\mathbf{A},\mathbf{i})$ close to $\left|{\mathbf{A}}\right|$ means that ${\mathbf{i}}$ has a high influence on ${\mathbf{A}}$, while an $CI(\mathbf{A},\mathbf{i})$ close to zero implies a negligible influence of ${\mathbf{i}}$ on ${\mathbf{A}}$. In other words, the higher the index, the higher the contagion effect of the event $\left\{X_{\mathbf i}>u\right\}$ on the region ${\mathbf{A}}$.

\begin{remark}{\rm{
 To gain some intuition for this measure, as a device for measuring dependence, consider two polar cases:

\begin{itemize}
	\item Case 1. If $X_{\mathbf {i}}$ is independent of $X_{\mathbf {j}}$, for each ${\mathbf {j}}\in {\mathbf {A}}$, then $CI(\mathbf{A},\mathbf{i})=0$.
	
	\item Case 2. If, for each ${\mathbf {j}}\in {\mathbf {A}}$, $X_{\mathbf {j}}$ and $X_{\mathbf {i}}$ are totally dependent, then $CI(\mathbf{A},\mathbf{i})=\left|{\mathbf{A}}\right|$.
\end{itemize}
}}
\end{remark}\vspace{0.3cm}

 \begin{remark}{\rm{
We can extend the $CI(\mathbf{A},\mathbf{i})$ to the contagion index from a region ${\mathbf{A}}$ to a region ${\mathbf{B}}$, as follows
\begin{equation*}
CI(\mathbf{A},\mathbf{B})=\lim_{u\uparrow 1} E\left(\sum_{\mathbf{j}\in \mathbf{A}}\indi_{\{F(X_{\mathbf j})>u\}}\left|\bigcup_{\mathbf{i}\in \mathbf{B}}\{F(X_{\mathbf i})>u\}\right.\right) .\label{SI}
\end{equation*}
This measure is related with the multivariate upper tail dependence coefficient (Schmidt (2002, \cite{schmidt}); Li (2009, \cite{li}) ; Ferreira (2011, \cite{fer})), defined as
$$
\lambda_{{\bf A},{\bf B}}=\lim_{u\uparrow 1}P\left(\bigcap_{{\bf j}\in {\bf A} }\{F(X_{\bf j})>u\}\left|\bigcap_{{\bf i}\in {\bf B} }\{F(X_{\bf i})>u\}\right.\right),
$$
in the following way
$$
CI(\mathbf{A},\mathbf{B})=\sum_{{\bf j}\in {\bf A}}\frac{\sum_{\emptyset\neq {\mathbf{J}}\subseteq {\mathbf{B}}}(-1)^{\left|{\mathbf{J}} \right|+1}\lambda_{{\mathbf{J}},\{{\mathbf{j}}\}}}{\epsilon_{\mathbf B}}.
$$
When we take ${\mathbf{A}}={\mathbf{B}}$, we obtain the fragility index (FI) of the region ${\mathbf{A}}$. The FI was introduced in Geluk \textit{et al}. (2007, \cite{geluk}) to measure the stability of a stochastic system. The system is called stable if $FI=1$, otherwise it is called fragile.
}}
\end{remark}

\section{Stability index of a region ${\mathbf{A}}$}

\pg In order to analyze the regional smoothness of a random field associated to a specific location we propose the following measure.

 \begin{definition}
 Let $X = \{X_{\bf{i}}\}_{{\bf{i}}\in\R^2}$ be a max-stable random field with unit Fréchet margins $F$ and $\mathbf{A}$ a region of $\mathbb{R}^2$. The stability index of the region $\mathbf{A}$ associated to a specific location ${\mathbf{i}}\in\mathbb{R}^2$, $SI({\mathbf{A}},{\mathbf{i}})$, is defined as
\begin{equation*}
SI(\mathbf{A},\mathbf{i})=\lim_{u\uparrow 1} E\left(\sum_{\mathbf{j}\in \mathbf{A}}\indi_{\{F(X_{\mathbf i})\leq u<F(X_{\mathbf j})\}}\left|\sum_{\mathbf{j}\in \mathbf{A}}\indi_{\{F(X_{\mathbf i})\leq u<F(X_{\mathbf j})\}}>0\right.\right). \label{SI}
\end{equation*}
\end{definition}\vspace{0.3cm}

The $SI({\mathbf{A}},{\mathbf{i}})$ is the conditional expected number of crossings (above a high threshold $u$) in ${\mathbf{A}}$ from a specific location ${\mathbf{i}}$, given that there is at least one crossing in ${\mathbf{A}}$ from ${\mathbf{i}}$.

If a max-stable random field ${\mathbf{X}}$ does not vary smoothly over a region ${\mathbf{A}}$, we will expect a large number of crossings of a high threshold in ${\mathbf{A}}$ associated to a specific location ${\mathbf{i}}$. A higher number of crossings signifies increased instability.

The next results highlight the connections between $SI(\mathbf{A},\mathbf{i})$ and the extremal coefficients.

\begin{proposition}
For any max-stable random field with unit Fréchet margins, ${\mathbf{i}}\in {\mathbb{R}}^2$ and  ${\mathbf{A}}\subset \mathbb{R}^2$, we have
$$
SI(\mathbf{A},\mathbf{i})=\frac{\sum_{\mathbf{j}\in \mathbf{A}}\epsilon_{\{{\mathbf {i}},{\mathbf{j}}\}}-\left|{\mathbf{A}}\right|}{\epsilon_{\{{\mathbf {i}}\}\cup{\mathbf{A}}}}=\frac{\left|{\mathbf{A}}\right|-CI(\mathbf{A},\mathbf{i})}{\epsilon_{\{{\mathbf {i}}\}\cup{\mathbf{A}}}}
$$\end{proposition}\vspace{0.3cm}

\bdem  Since
$$
E\left(\sum_{\mathbf{j}\in \mathbf{A}}\indi_{\{F(X_{\mathbf i})\leq u<F(X_{\mathbf j})\}}\left|\sum_{\mathbf{j}\in \mathbf{A}}\indi_{\{F(X_{\mathbf i})\leq u<F(X_{\mathbf j})\}}>0\right.\right)=\frac{\sum_{\mathbf{j}\in \mathbf{A}}P(F(X_{\mathbf i})\leq u<F(X_{\mathbf j}))}{1-P(F(X_{\mathbf i})\leq u,\bigcap_{\mathbf{j}\in \mathbf{A}} \{F(X_{\mathbf j})\leq u\})},
$$
it follows that
\begin{eqnarray*}
SI(\mathbf{A},\mathbf{i})&=&\lim_{u\uparrow 1} \frac{\left|\mathbf{A}\right|u-\sum_{\mathbf{j}\in \mathbf{A}}u^{\epsilon_{\{\mathbf{i},\mathbf{j}\}}}}{1-u^{\epsilon_{\{\mathbf{i}\}\cup \mathbf{A}}}}=\frac{\sum_{\mathbf{j}\in \mathbf{A}}\epsilon_{\{\mathbf{i},\mathbf{j}\}}-\left|\mathbf{A}\right|}{\epsilon_{\{{\mathbf {i}}\}\cup{\mathbf{A}}}}.
\end{eqnarray*}
\edem\vspace{0.3cm}

\begin{proposition}
Under the conditions of Proposition 3.1., we have
$$
\frac{\sum_{\mathbf{j}\in \mathbf{A}}\epsilon_{\{\mathbf{i},\mathbf{j}\}}-\left|\mathbf{A}\right|}{\left|\mathbf{A}\right|+1}\leq SI(\mathbf{A},\mathbf{i})\leq\frac{\sum_{\mathbf{j}\in \mathbf{A}}\epsilon_{\{{\mathbf {i}},{\mathbf{j}}\}}-\left|{\mathbf{A}}\right|}{\bigvee_{{\mathbf {j}}\in {\mathbf {A}}}\epsilon_{\{{\mathbf {i}},{\mathbf{j}}\}}}.
$$\end{proposition}\vspace{0.3cm}

\bdem  Just observe that $\epsilon_{\{\mathbf {i}\}\cup{\mathbf{A}}}\leq\left|\mathbf{A}\right|+1$ and $\epsilon_{\{{\mathbf {i}}\}\cup{\mathbf{A}}}\geq \bigvee_{{\mathbf {j}}\in {\mathbf {A}}}\epsilon_{\{{\mathbf {i}},{\mathbf{j}}\}}$.
\edem\vspace{0.3cm}


\bigskip


\section{An M4 random field}


In this section we derive the contagion and stability indices for a particular class of max-stable processes known as multivariate maxima of moving maxima, or M4 processes for short (Smith and Weissman (1996), \cite{smi2}), which are particularly well applicable in a time series or spatial process context. An M4 random field ${\bf X}=\left\{X_{\mathbf{i}}\right\}_{{\mathbf{i}}\in \mathbb Z^2}$ is defined as
\begin{equation}
X_{{\mathbf{i}}}=\max_{l\geq 1}\max_{-\infty<m<+\infty} a_{lm{\mathbf{i}}}Z_{l,1-m},\quad {\mathbf{i}}\in \mathbb{Z}^2,\label{campo_M4}
\end{equation}
where $\{Z_{l,n}\}_{l\geq 1, n\in \mathbb{Z}}$ is a family of independent unit Fréchet random variables and, for each ${\mathbf{i}}\in \mathbb{Z}^2$, $\{a_{lm{\mathbf{i}}}\}_{l\geq 1, m\in \mathbb{Z}}$ are non-negative constants such that $\displaystyle{\sum_{l=1}^{+\infty}\sum_{m=-\infty}^{+\infty}a_{lm{\mathbf{i}}}=1}.$
By considering that  the distribution of $(X_{{\mathbf{i}}_1},\ldots,X_{{\mathbf{i}}_k})$ is characterized by the copula
\begin{eqnarray}
C(u_{{\mathbf{i}}_1},\ldots,u_{{\mathbf{i}}_k})=\prod_{l=1}^{+\infty}\prod_{m=-\infty}^{+\infty}\bigwedge_{{\mathbf{i}}\in \{{\mathbf{i}}_1,\ldots,{\mathbf{i}}_k\}}u_{{\mathbf{i}}}^{a_{lm{\mathbf{i}}}},\quad u_{{\mathbf{i}}_j}\in [0,1],\ j=1,\ldots,k,\label{copula}
\end{eqnarray}
it was shown in Fonseca \textit{et al.} (2012, \cite{cec}) that the random field ${\bf X}=\left\{X_{\mathbf{i}}\right\}_{{\mathbf{i}}\in \mathbb Z^2}$ is max-stable and the exponent function of the distribution of $(X_{{\mathbf{i}}_1},\ldots,X_{{\mathbf{i}}_k})$ is given by
\begin{eqnarray*}
V_{{\bf{A}}}(x_1,\ldots,x_k)=\sum_{l=1}^{+\infty}\sum_{m=-\infty}^{+\infty} \bigvee_{j=1}^{k}\left(x_j^{-1}a_{lm\mathbf{i}_j}\right),\ x_j\in \R,\ j=1,\ldots,k,\ {\bf{A}}=\{{\bf{i}}_1,\ldots,{\bf{i}}_k\}.
\end{eqnarray*}
So
$$
CI({\mathbf{A}},{\mathbf{i}})=2\left|{\mathbf{A}}\right|-\sum_{{\mathbf{j}}\in {\mathbf{A}}}\sum_{l=1}^{+\infty}\sum_{m=-\infty}^{+\infty}(a_{lm{\mathbf{i}}}\vee a_{lm{\mathbf{j}}})
$$
and
$$
SI({\mathbf{A}},{\mathbf{i}})=\frac{\sum_{{\mathbf{j}}\in {\mathbf{A}}}\sum_{l=1}^{+\infty}\sum_{m=-\infty}^{+\infty}(a_{lm{\mathbf{i}}}\vee a_{lm{\mathbf{j}}})-\left|{\mathbf{A}}\right|}{\sum_{l=1}^{+\infty}\sum_{m=-\infty}^{+\infty}\left(\bigvee_{{\mathbf{j}}\in {\mathbf{A}}}a_{lm{\mathbf{j}}}\vee a_{lm{\mathbf{i}}}\right)}.
$$

\bigskip

To illustrate the computation of the contagion and stability indices we shall consider, in what \linebreak follows, examples with a finite number of signature patterns $(1\leq l\leq L)$ and a finite range of sequential dependencies $(M_1\leq m\leq M_2)$.

\begin{ex}

{\rm{Lets consider that for each location ${\mathbf{i}}\in \mathbb{Z}^2$ with even abscissa we have \linebreak $a_{11{\mathbf{i}}}=\frac{4}{5}$, $a_{12{\mathbf{i}}}=\frac{1}{5}$ and otherwise $a_{11{\mathbf{i}}}=\frac{1}{4}=1-a_{12{\mathbf{i}}}.$ The values of
$(a_{11{\mathbf{i}}},a_{12{\mathbf{i}}})$ determine the \linebreak moving pattern or signature pattern  of the random field, which in this case corresponds to one pattern ($L=1$).

\begin{figure}[!htb]
\begin{center}
\includegraphics[scale=0.45]{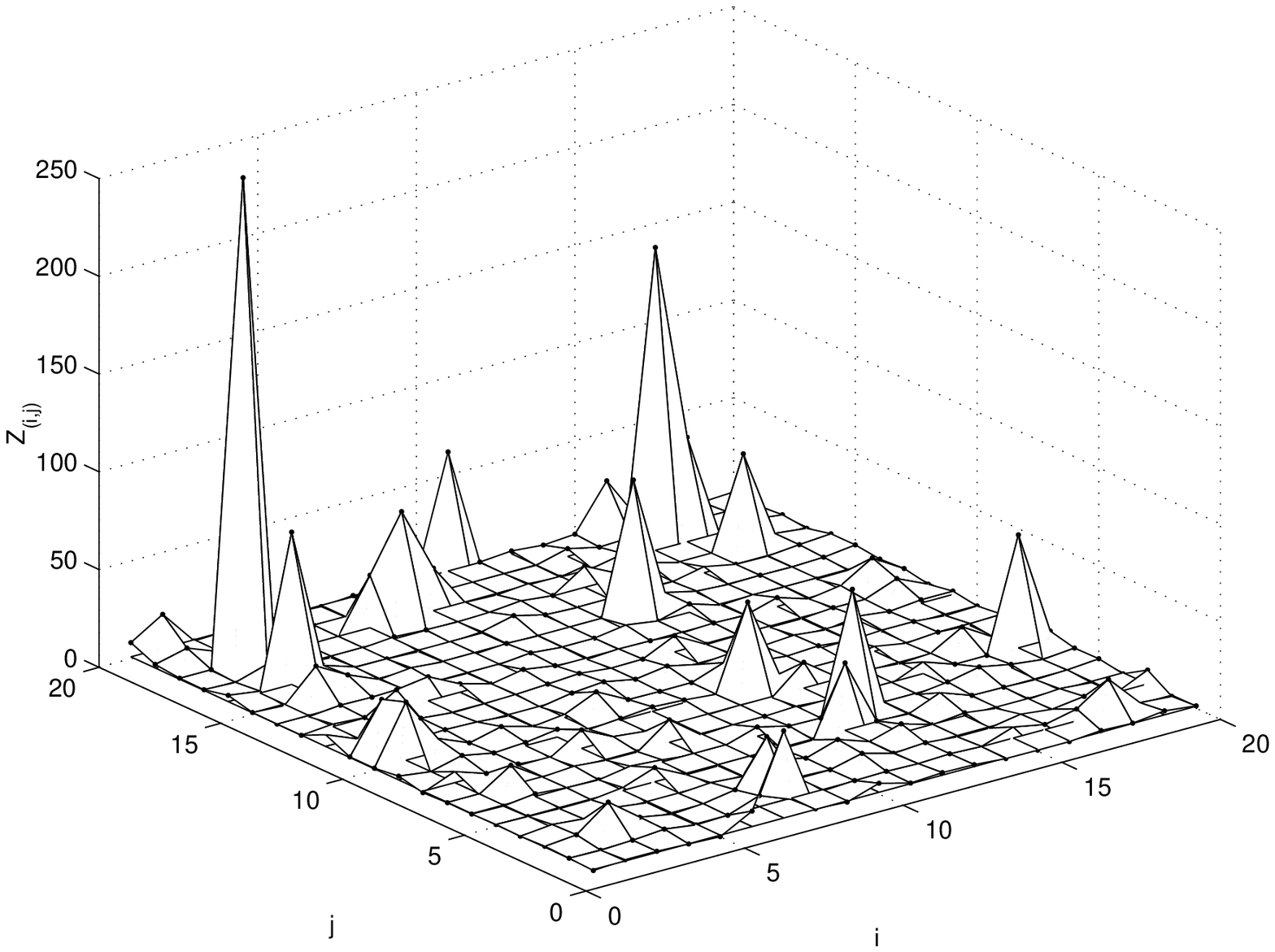}
\includegraphics[scale=0.45]{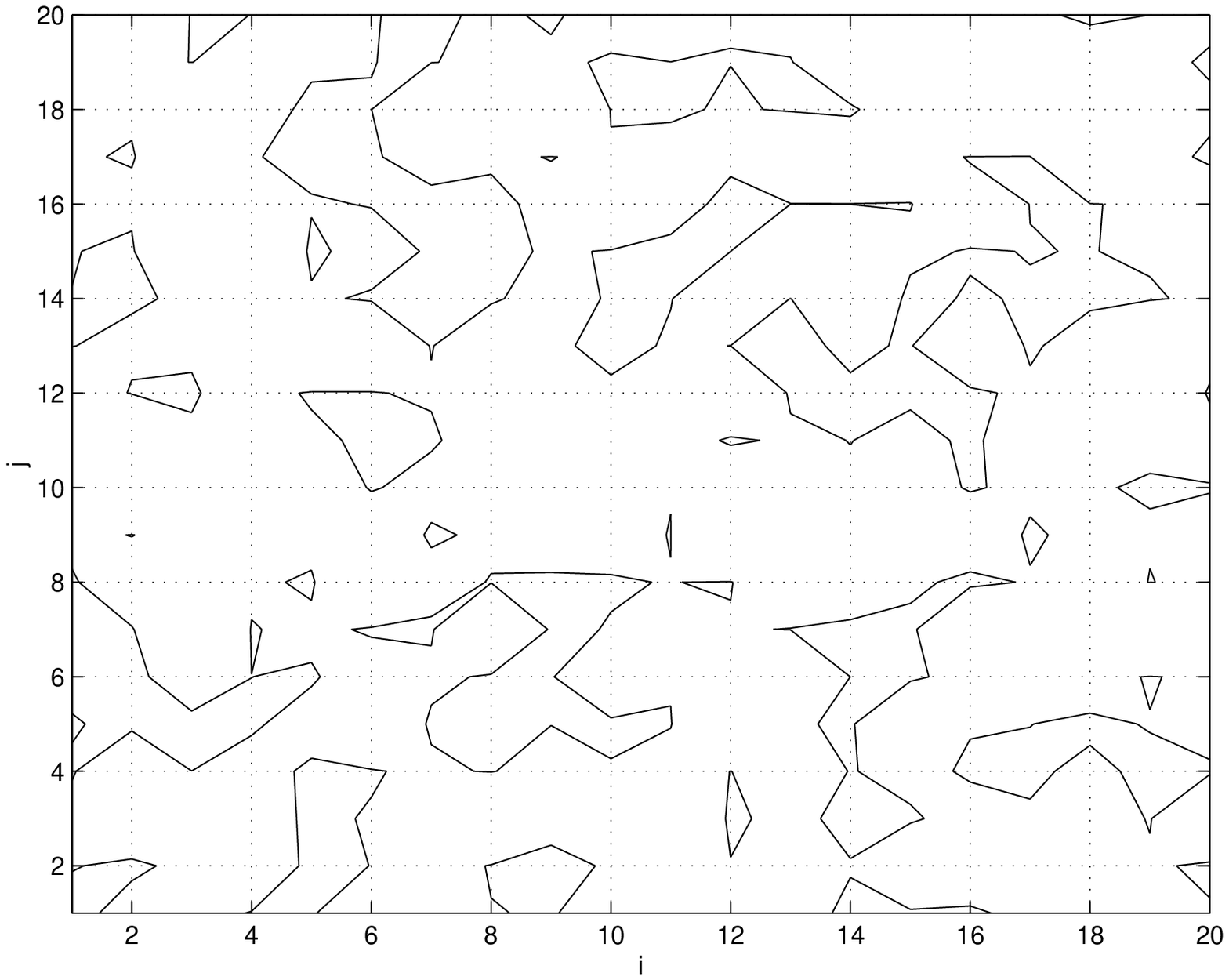}
\vspace{-0.5cm}\caption{\small Simulation of the M4 random field with $L=1$, $1\leq m\leq 2$ and for constants defined in Example 4.1.}
\end{center}
\label{fig:figura1}
\end{figure}

Let ${\mathbf{A}}=\{s_{j}({\mathbf{i}}):{\mathbf{i}}=(3,3)\wedge j\in\{1,2,\ldots,8\}\}$, where $s_{j}({\mathbf{i}})$, $j=1,\ldots,8$, denote the neighbors of ${\mathbf{i}}$ defined as follows:
\begin{eqnarray*}
s_1({\mathbf{i}})=(i_1+1,i_2) \ \ s_2({\mathbf{i}})={\mathbf{i}}+{\mathbf{1}} \ \ s_3({\mathbf{i}})=(i_1,i_2+1) \ \ s_4({\mathbf{i}})=(i_1-1,i_2+1)\\
s_5({\mathbf{i}})=(i_1-1,i_2) \ \ s_6({\mathbf{i}})={\mathbf{i}}-{\mathbf{1}} \ \ s_7({\mathbf{i}})=(i_1,i_2-1) \ \ s_8({\mathbf{i}})=(i_1+1,i_2-1).
\end{eqnarray*}
The matrix of the bivariate extremal coefficients, $\epsilon_{\{s_j({\mathbf{i}}),{\mathbf{i}}\}},$ $j=1,\ldots,8$, provide insight into the likelihood of contagion from ${\mathbf{i}}=(3,3)$ to its neighbors $s_{j}({\mathbf{i}})$, $j=1,\ldots,8$, although without specifying the size of the impact which is given by $CI({\mathbf{A}},{\mathbf{i}})$. We obtain
\begin{equation*}
\left[
                         \begin{array}{ccc}
                           \epsilon_{\{s_4(\mathbf{i}),\mathbf{i}\}}  & \epsilon_{\{s_3(\mathbf{i}),\mathbf{i}\}} & \epsilon_{\{s_2(\mathbf{i}),\mathbf{i}\}} \\
                           \epsilon_{\{s_5(\mathbf{i}),\mathbf{i}\}} & \epsilon_{\{\mathbf{i},\mathbf{i}\}} & \epsilon_{\{s_1(\mathbf{i}),\mathbf{i}\}} \\
                           \epsilon_{\{s_6(\mathbf{i}),\mathbf{i}\}} & \epsilon_{\{s_7(\mathbf{i}),\mathbf{i}\}} & \epsilon_{\{s_8(\mathbf{i}),\mathbf{i}\}} \\
                         \end{array}
                       \right]= \left[
                         \begin{array}{ccc}
                           \frac{31}{20}  & 1 & \frac{31}{20} \\
                           \frac{31}{20} & 1 & \frac{31}{20} \\
                           \frac{31}{20} & 1 & \frac{31}{20} \\
                         \end{array}
                       \right]
\end{equation*} \vspace{0.3cm}
and
$$
CI(\mathbf{A},\mathbf{i})=4,7.
$$
The stability index of the region ${\mathbf{A}}$ associated with $\mathbf{i}=(3,3)$ is given by
$$
SI(\mathbf{A},\mathbf{i})=\frac{66}{31}.
$$
}}
\end{ex}

\begin{ex}
{\rm{Now, we shall consider one example with two signature patterns ($L=2$).

Lets assume that for each location ${\mathbf{i}}\in \mathbb{Z}^2$ we have $a_{11{\mathbf{i}}}=a_{12{\mathbf{i}}}=a_{13{\mathbf{i}}}=\frac{1}{5},$ \linebreak $a_{21{\mathbf{i}}}=a_{22{\mathbf{i}}}=\frac{1}{10}$, $a_{23{\mathbf{i}}}=\frac{1}{5}$ if both coordinates are odd and $a_{11{\mathbf{i}}}=\frac{1}{4}, a_{12{\mathbf{i}}}=a_{13{\mathbf{i}}}=\frac{1}{8},$ \linebreak $a_{21{\mathbf{i}}}=a_{22{\mathbf{i}}}=a_{23{\mathbf{i}}}=\frac{1}{6}$ otherwise. The values of $(a_{11{\mathbf{i}}},a_{12{\mathbf{i}}},a_{13{\mathbf{i}}})$ and $(a_{21{\mathbf{i}}},a_{22{\mathbf{i}}},a_{23{\mathbf{i}}})$ define the two signature patterns of the random field.

\begin{figure}[!htb]
\begin{center}
\includegraphics[scale=0.45]{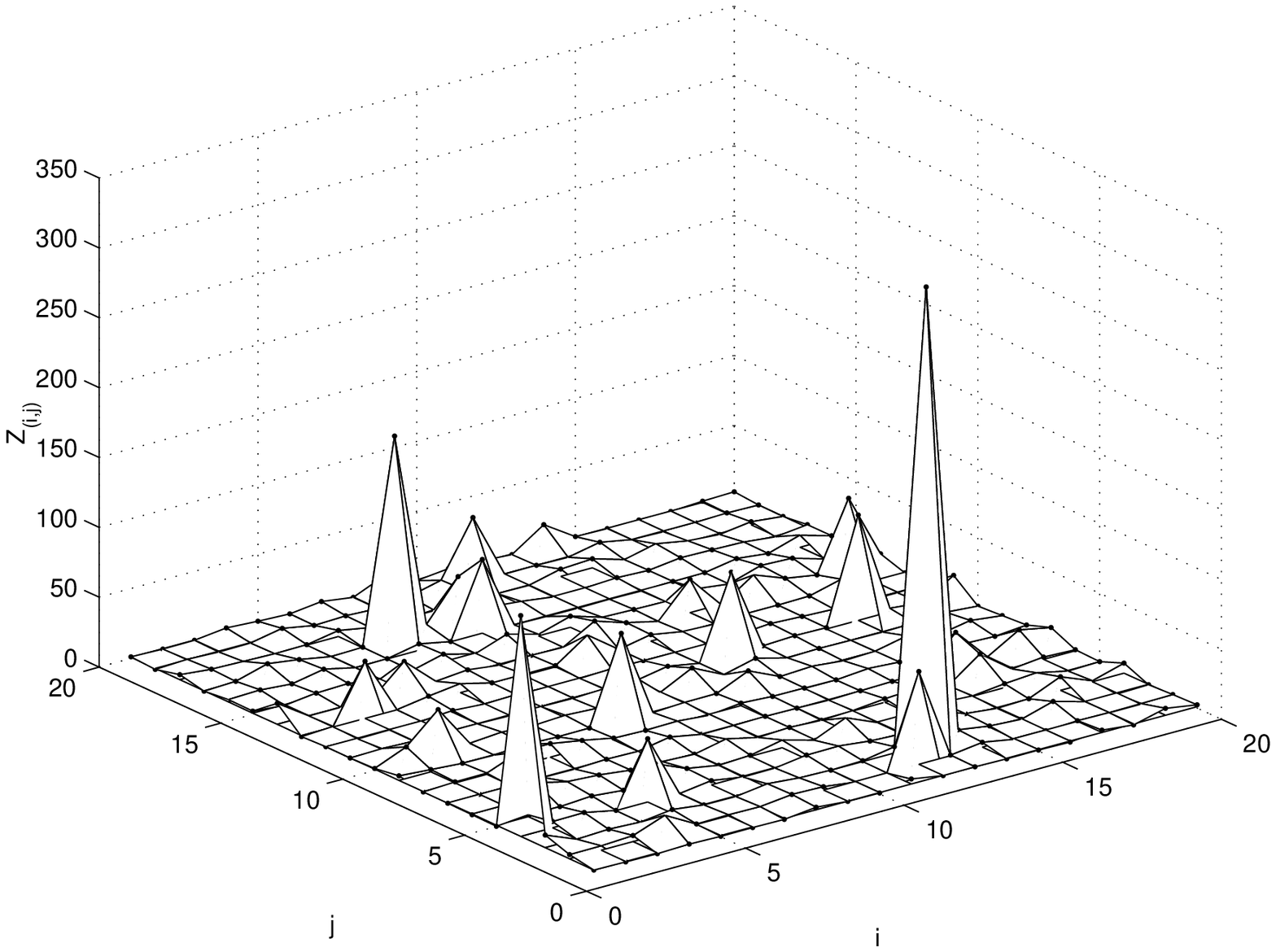}
\includegraphics[scale=0.45]{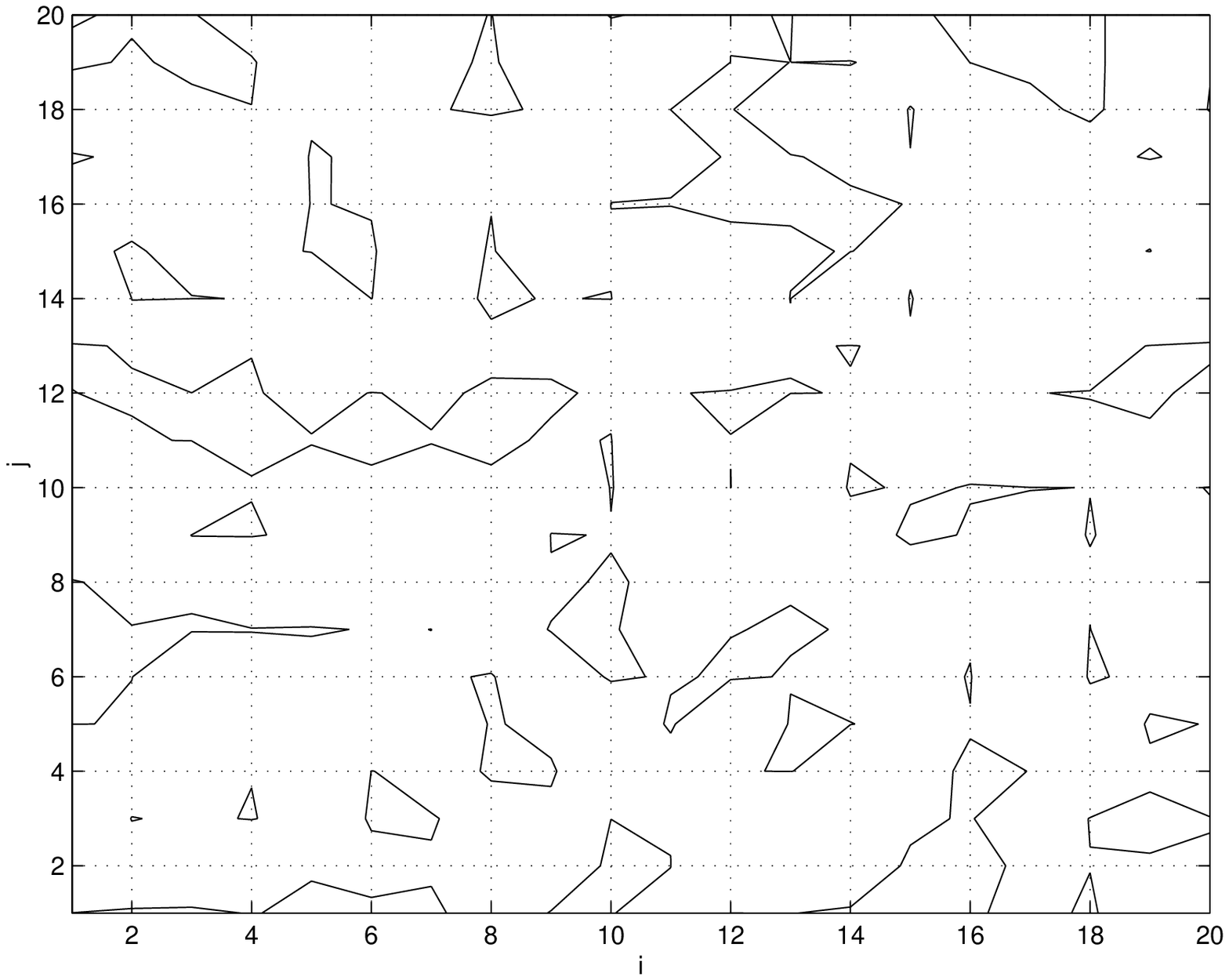}
\vspace{-0.5cm}\caption{\small Simulation of the M4 random field with $L=2$, $1\leq m\leq 2$ and for constants defined in Example 4.2.}
\end{center}
\label{fig:figura2}
\end{figure}

For ${\mathbf{A}}=\{(2,4),(3,4), (4,4), (5,4)\}$ and ${\mathbf{i}}=(3,3)$ we obtain $CI({\mathbf{A}},{\mathbf{i}})=\frac{49}{15}$ and $SI({\mathbf{A}},{\mathbf{i}})=\frac{44}{71}$.
}}
\end{ex}

We next focus on the estimation of the stability and contagion indices.

\section{Estimation}

As previously stated, the contagion and stability indices relate with the extremal coefficients of Tiago de Oliveira (1962/62, \cite{tiago}) and Schlather and Tawn (2003, \cite{sch}), which can be expressed by the exponent function given in (1).

There are several references in literature on the estimation of the exponent function. For a survey we refer to Krajina (2010, \cite{kraj}) and Beirlant (2004, \cite{beir}).

It is known that parametric estimation methods are efficient if the distribution model under consideration is true, but they suffer from biased estimates otherwise. Non parametric estimation procedures avoid this type of model error. However, they are usually based on an arbitrarily chosen parameter $k$ corresponding to the number of top order statistics to be used on the estimation of a high quantile of $F$, which relates to the usual variance-bias problem: if $k$ is too small, then the estimator tends to have a large variance, where if $k$ is too large, then the bias tends to dominate. Some methods of choosing an optimal $k$ are discussed in Einmahl \textit{et al.} (2006, \cite{einmahl}).

In order to overcome the problem of the optimal choice of $k$, Ferreira and Ferreira (2011, \cite{fer2}) developed another approach. Based on the following relation
$$
\epsilon_{\bf{A}}=V_{\bf{A}}(1,1,\ldots,1)=\frac{E(M({\bf A}))}{1-E(M({\bf A}))}, \ \ \ {\text {where}} \ \ \ M({\bf A})=\bigvee_{{\bf {i}}\in {\bf {A}}}F_{\bf i}(X_{\bf i}),
$$
the estimator of $\epsilon_{\bf{A}}$ proposed in Ferreira and Ferreira (2011, \cite{fer2}) is defined as
$$
\widehat{\epsilon}_{\bf{A}}=\frac{\overline{M({\bf{A}})}}{1-\overline{M({\bf{A}})}},
$$
where $\overline{M({\bf{A}})}$ is the sample mean,
$$
\overline{M({\bf{A}})}=\frac{1}{n}\sum_{j=1}^n \bigvee_{{\bf {i}}\in {\bf {A}}}\widehat{F_{\bf i}}(X_{\bf i}^{(j)})
$$
and $\widehat{F_{\bf i}}$, ${\bf {i}}\in {\bf {A}}$, is the (modified) empirical distribution function of $F_{\bf i}$,
$$
\widehat{F_{\bf i}}(u)=\frac{1}{n+1}\sum_{j=1}^n\indi_{\left\{X_{\bf i}^{(j)}\leq u\right\}}.
$$

With this estimator of the extremal coefficient and the relations established in Propositions 2.1 and 3.1 we propose, respectively, the following estimators for the contagion index $CI({\bf A},{\bf i})$ and the stability index $SI({\bf A},{\bf i})$,
$$
\widehat{CI}({\bf A},{\bf i})=2\left|{\bf A}\right|-\sum_{{\bf j}\in {\bf A}}\widehat{\epsilon}_{\left\{{\bf i},{\bf j}\right\}}
$$
and

$$
\widehat{SI}({\bf A},{\bf i})=\frac{\sum_{{\bf j}\in {\bf A}}\widehat{\epsilon}_{\left\{{\bf i},{\bf j}\right\}}-\left|{\bf A}\right|}{\widehat{\epsilon}_{\left\{{\bf i}\right\}\cup{\bf A}}},
$$
which are consistent given the consistency of the estimators $\widehat{\epsilon}_{\left\{{\bf i},{\bf j}\right\}}$ and $\widehat{\epsilon}_{\left\{{\bf i}\right\}\cup{\bf A}}$ already stated in Ferreira and Ferreira (2011, \cite{fer2}).

The results of the application of these estimators to the Examples 4.1 and 4.2 are presented in the following tables:

 \begin{table}[h]
\begin{center}
\begin{tabular}{|c|c|c|c|}  \cline{2-4}
\multicolumn{1}{c|}{}  & {  \emph{CI}}  &  {  $\widehat{CI}$ } & {  \emph{MSE} } \\ \hline
Example 4.1  & 4.7 & 4.70806  &  0.01259 \\ \hline
Example 4.2 & $3.2667$ & 3.26684  & 0.0006 \\\hline
\end{tabular}
\caption{Results with 100 replications of 1000 i.i.d. max-stable M4 random fields of the Examples 4.1 and 4.2 where $CI$ denotes the true values of the contagion index, $\widehat{CI}$ the estimated values and $MSE$ the estimated mean squared error.} \label{tab:IC}
\end{center}
\end{table}

\begin{table}[h]
\begin{center}
\begin{tabular}{|c|c|c|c|} \cline{2-4}
  \multicolumn{1}{c|}{} & { \emph{SI}} & { $\widehat{SI}$} & { \emph{MSE}} \\ \hline
Example 4.1    & $2.12903$ & 2.1222 & 0.00234  \\ \hline
Example 4.2 & $0.61972$ & 0.62029 & 0.00034 \\\hline
\end{tabular}
\caption{Results with 100 replications of 1000 i.i.d. max-stable M4 random fields of the Examples 4.1 and 4.2 where $SI$ denotes the true values of the stability index, $\widehat{SI}$ the estimated values and $MSE$ the estimated mean squared error.} \label{tab:SI}
\end{center}
\end{table}

 As we can see the estimated values are very close to the true values of the coefficients. These results show that the simple non-parametric estimators $\widehat{CI}$ and $\widehat{SI}$ are a promising tool for assessing regional contagion effects and regional smoothness for these random fields.

\section{An application to precipitation data}

\pg We now illustrate the estimation of the contagion and stability indices through an application to annual maxima values of daily maxima precipitation recorded over 32 years, in six Portuguese stations (Figure 3), provided by the Portuguese National System of Water Resources (http://snirh.pt).

\begin{figure}[!htb]
\begin{center}
\includegraphics[scale=0.25]{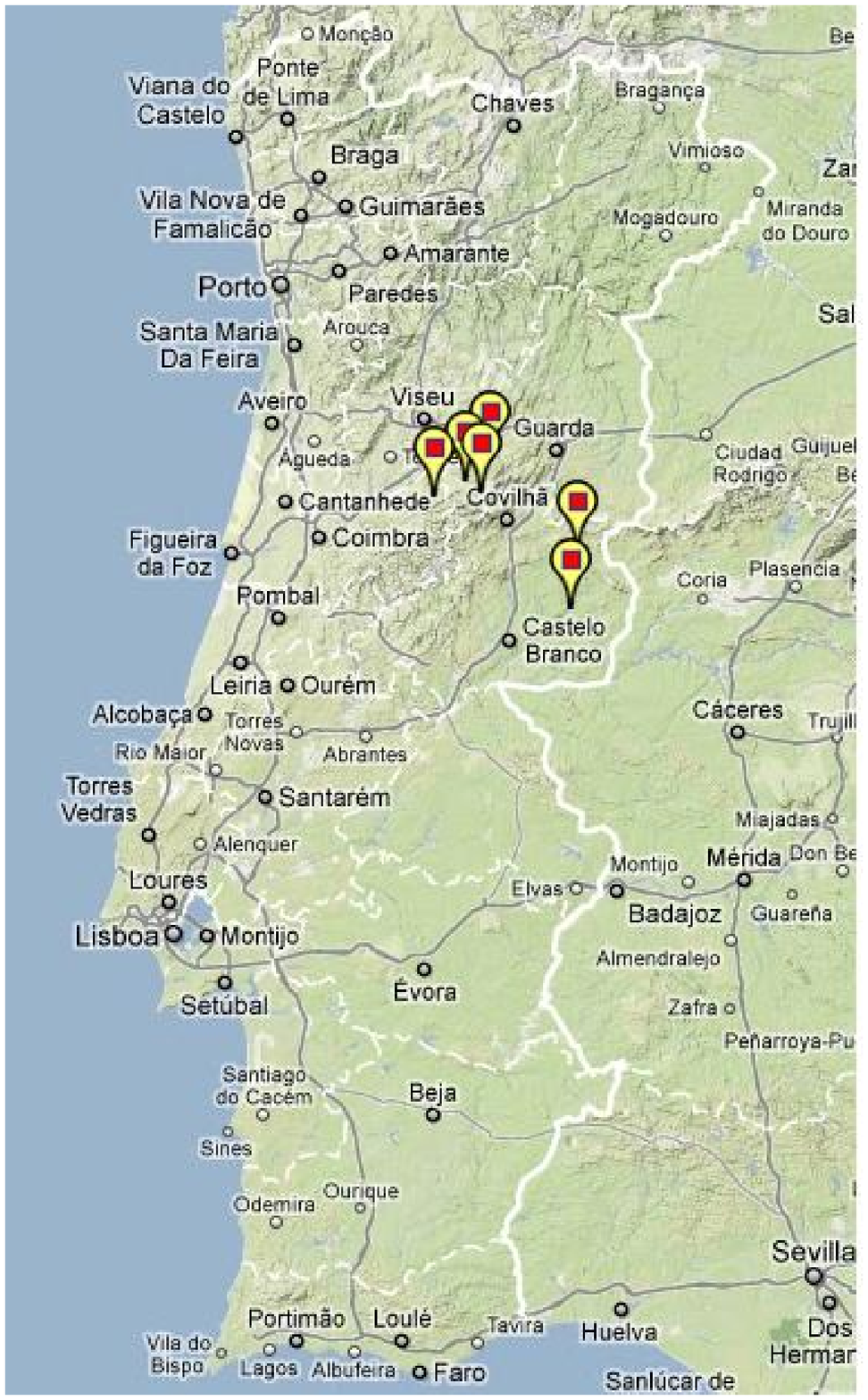} \ \ \ \ \ \
\includegraphics[scale=0.35]{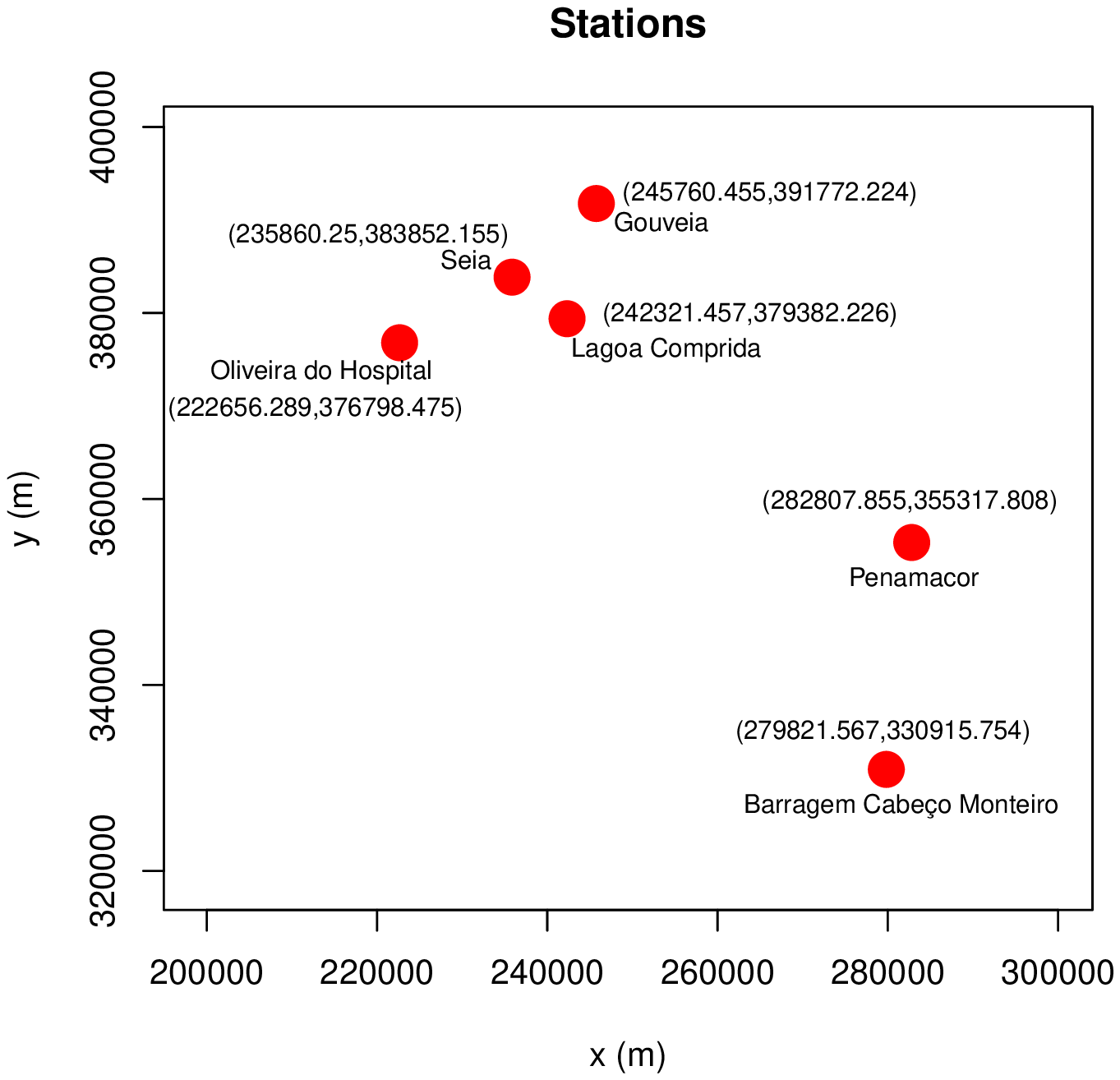}
\vspace{-0.5cm}\caption{\small The locations of the stations where precipitation data were collected, obtained from Portuguese National System of Water Resources (left) and their representation in Lambert coordinates (right).}
\end{center}
\label{fig:figura3}
\end{figure}

Since the data are maxima over a long period of time, we assumed that they are independent over the years in each location. We also assumed that the random field is max-stable with unknown marginal distributions so data were previously transformed at each site so that they have a standard Fréchet distribution.

The estimated values of the contagion and stability indices from Lagoa Comprida (located at the highest altitude) to the regions
 \{Gouveia, Oliveira do Hospital, Seia\} and \{Penamacor, Barragem Cabeço Monteiro\} are presented in Table 3.

\begin{table}[h]
\begin{center}
\begin{tabular}{|c||c|c|} \hline
${\bf A}$ & {$\widehat{CI}({\bf A}, {\text {Lagoa \ Comprida}})$} & {$\widehat{SI}({\bf A}, {\text {Lagoa \ Comprida}})$}  \\ \hline \hline
$\{{\text{Gouveia, Oliveira \ do \ Hospital, Seia}}\}$    & 0.96688 & 0.8959   \\ \hline
$\{{\text{Penamacor, Barragem \ Cabeço \ Monteiro}}\}$ & 0.00887 & 0.75741  \\\hline
\end{tabular}
\caption{Estimates of the contagion and stability indices from Lagoa Comprida to the regions \{Gouveia, Oliveira do Hospital, Seia\} and \{Penamacor, Barragem Cabeço Monteiro\}.} \label{tab:3}
\end{center}
\end{table}

The results suggest that Lagoa Comprida has a higher influence on the region \{Gouveia, Oliveira do Hospital, Seia\} in terms of precipitation amounts and this region is smoother when compared to region \{Penamacor, Barragem Cabeço Monteiro\}.



\end{document}